\documentclass[12pt,leqno]{amsart}
\usepackage{amscd}
\usepackage{amssymb}
\usepackage{amsfonts}
\usepackage{latexsym}
\usepackage{verbatim}
\theoremstyle{plain}
\DeclareMathOperator{\End}{End}
\newtheorem{theorem}{Theorem}[section]

\newtheorem{lemma}[theorem]{Lemma}

\newtheorem{cor}[theorem]{Corollary}

\newcommand{\ov}[1]{\overline{#1}}

\newcommand{\ka}[1]{\kappa{#1}}
\newcommand{\de}[2]{\frac{\partial #1}{\partial #2}}
\newcommand{\n}{\nabla}
\begin{document}
\title[A generalization of the  normal holomorphic frames]{A generalization of the normal holomorphic frames in
Symplectic Manifolds}
\author{Luigi Vezzoni}
\date{\today}
\address{Dipartimento di Matematica ''L. Tonelli''\\ Universit\`a di Pisa\\
Via Buonarroti 2\\
56127 Pisa\\ Italy} \email{vezzoni@mail.dm.unipi.it}
\thanks{This work was supported by G.N.S.A.G.A.
of I.N.d.A.M.}
\begin{abstract}In  this paper
we give a generalization of the normal
holomorphic frames in the symplectic manifolds and find conditions for
the integrability of complex structures.
\end{abstract}
\maketitle
\section{Introduction}
One of the possible ways of studying symplectic manifolds is to fix a
complex structure calibrated by the symplectic form and then to
consider a structure which is the natural generalization of
K\"ahler structure (sometimes in
literature this structure is called  almost-K\"ahler
structure). \\
Moreover we can try to
extend K\"ahler features in symplectic manifolds and find
conditions on the Riemannian invariants (like the curvature e.g.) which
force  the integrability of the complex structure. These problems are
object of research of many papers (see e.g. \cite{Ab}, \cite{I},
\cite{L}, \cite{SEk}). \\
Starting from an idea of P. de Bartolomeis and  A. Tomassini (see
\cite{J}), in this paper we generalize the normal holomorphic
frames, characteristic of the K\"ahler manifolds, to the symplectic
case (see theorem 1); we call these frames
\sl{generalized normal holomorphic frames}\rm.\\
As  application of the existence of generalized normal holomorphic
frames we prove that if the (0,1)-part of the covariant derivative of
the (1,2)-tensor
$$
B(X,Y)=J(\n_XJ)Y-(\n_{JX}J)Y
$$
vanishes, then the complex structure $J$ is integrable (see theorem 2).\\

The paper is organized as follows.
After some preliminaries, in  section 3 we give the proof of the
existence of generalized normal holomorphic frames and in section
4 we prove theorem 2. Finally we give another sufficient condition for
the integrability of $J$. Namely, we show that if the complex
structure satisfies
$$
(\n J)\n=0\,,$$
then $J$ is holomorphic.
\section{Preliminaries}
Let $M$ be a $2n$-dimensional manifold: a \emph{complex structure} on $M$ is
a smooth  section $J$ of the bundle $\End(TM)$ such that $J^2=-Id$.
The couple $(M,J)$ is
called a \emph{complex manifold}.\\
A complex structure gives a natural split of $TM\otimes
\mathbb{C}$ in $TM^{(1,0)}\oplus TM^{(0,1)}$, where $TM^{(1,0)}$ and
$TM^{(0,1)}$
are the eigenspaces relatively to $i$ and $-i$.\\
A Riemannian $J$-invariant metric $g$ on
a complex manifold $(M,J)$ is said to be Hermitian and the triple
$(M,g,J)$ is by definition
an
Hermitian manifold. An Hermitian metric $g$
induces a non-degenerate 2-form $\kappa$ on $M$, given by
$$
\kappa(J\cdot,\cdot)=g(\cdot,\cdot)\,.
$$
It is well known the fundamental Hermitian relation:
\begin{equation}
2g((\nabla_{X}J)Y,Z)=d\ka(X,JY,JZ)-d\ka(X,Y,Z)+
g(N_{J}(Y,Z),JX)
\end{equation}
(see e.g. \cite{F}),
where $\n$ is the Levi-Civita connection of $g$  and
$N_J$ is the Nijenhuis tensor of $J$, i.e. the (1,2)-tensor defined on $M$
by
$$
N_J(X,Y):=[JX,JY]-J[JX,Y]-J[X,JY]-[X,Y]\,.
$$

A holomorphic manifold $M$ has a natural complex structure.
A complex structure $J$ is
said to be integrable if it is induced by a holomorphic structure.
In view of
Newlander-Nirenberg theorem a complex structure is
integrable if and only if $N_J$ vanishes.\\

A \emph{symplectic manifold} is a pair $(M,\kappa)$, where $\kappa$ is a
non-degenerate closed 2-form (i.e. $\kappa^n\neq 0,\;\; d\kappa=0$) .\\
A complex
structure $J$ on a symplectic manifold  $(M,\kappa)$
is said to be $\kappa$-calibrated if
\begin{enumerate}
\item[1.] $\kappa(J\cdot,J\cdot)=\kappa(\cdot,\cdot)$ ($\kappa$ is $J$-invariant);
\vspace{0.1 cm}
\item[2.] $g_J(\cdot,\cdot):=\kappa(J\cdot,\cdot)$ is a positive defined tensor.
\end{enumerate}
Obviously $g_J$ is a Hermitian metric on $M$.
Let us denote by $\mathcal{C}_{\kappa}(M)$ the set of the $\kappa$-calibrated
complex structures on $M$: it is well known that
$\mathcal{C}_{\kappa}(M)$ is always a non-empty  set (see e.g. \cite{A}).\\

 A \emph{K\"ahler structure} on $M$ is a blend of three components
$(g,J,\kappa)$, where $\kappa$ is a symplectic structure,
$J$ is an integrable $\kappa$-calibrated complex structure and
$g=g_J$ (for the general theory about  K\"ahler manifolds we remand to
\cite{GFE}, \cite{G-H}, \cite{F} and \cite{H}).\\
It's known that a holomorphic Hermitian manifold $(M,g,J)$
is K\"ahler if and
only if around every point $o$ in $M$ there exists a normal
holomorphic coordinate system. Namely, if $o$ is any  point of $M$,
there exists a coordinate system
$\{z_1,\dots,z_n\}$  around $o$ such that, if
$G_{ij}:=g(\partial/\partial z_i,
\partial/\partial{\ov{z}_j})$, then
\begin{enumerate}
\item[1.] $G_{ij}(o)=\delta_{ij}$ ;
\vspace{0.1 cm}
\item[2.] $d G_{ij}[o]=0$ .
\end{enumerate}
A system of normal coordinates induces a (1,0)-frame $\{Z_1,\dots,Z_n\}$
around $o$ such that
$$
\n_{Z_i}Z_k(o)=0\,,\qquad \n_{\ov{Z}_i}Z_k(o)=0\,,
$$
for $1\leq i,k\leq n$, where $\n$ is the Levi-Civita connection of $g$.
We call such frames \emph{ ``normal holomorphic frames''}.

The aim of this paper is to construct a generalization of the normal
holomorphic frames in the symplectic manifolds and find some
non-obvious conditions which imply $N_J=0$.\\
\newline
At first we note that, in the  symplectic case, (1) reduces
to
\begin{equation}
2g((\nabla_{X}J)Y,Z)=g(N_{J}(Y,Z),JX)\,;
\end{equation}
so we have:
\begin{cor}
Let $(M,\kappa)$ be a symplectic manifold, $J\in\mathcal{C}_{\kappa}(M)$,
$g=g_J$ and $\n$ the Levi-Civita  connection of $g$.\\
Then
$\nabla_{\ov{Z}_1}Z_2$ is a (1,0)-vector field on $M$ for any
$Z_1,Z_2$ of  type (1,0).
\end{cor}
\begin{proof}
It's easy to check that:
\begin{equation}
N_J(Z_1,\ov{Z}_2)=0\,,
\end{equation}
for any $Z_1,Z_2$ $(1,0)$-vector fields.\\
Therefore, from the complex extension of (2), we have:
\begin{equation*}
\begin{aligned}
&g((\n_{\ov{Z}_1}J)Z_2,Z_3)=0\,,\\
&g((\n_{\ov{Z}_1}J)Z_2,\ov{Z}_3)=0.
\end{aligned}
\end{equation*}
Hence
$$
0=(\n_{\ov{Z}_1}J)Z_2=\n_{\ov{Z}_1}JZ_2-J\n_{\ov{Z}_1}Z_2=
i\n_{\ov{Z}_1}Z_2-J\n_{\ov{Z}_1}Z_2\,.
$$
\end{proof}
Now we have an obstruction to generalize normal holomorphic frames :
\begin{cor}
Let $(M,\ka)$ be a symplectic manifold,
$J\in\mathcal{C}_{\kappa}(M) $ and $g=g_J$.\\
Assume that for every point $\;o\in M$ there
exists a local complex (1,0)-frame \\$\{Z_{1},\dots,Z_{n}\}$,
around $o$, such that:
\begin{enumerate}
\item[1.] $\nabla_{Z_i}\ov{Z}_k(o)=0$, $1\leq i,k \leq n$;
\vspace{0.1 cm}
\item[2.] The real underlying frame $\{X_1,\dots,X_n\}$ satisfies
$$
\n_{X_{i}}X_{k}(o)=0\,,\,\n_{JX_{i}}X_k(o)=0\,,\,1 \leq i,k \leq n\,,
$$
\end{enumerate}
then $(g,J,\kappa)$ is a K\"ahler structure on $M$.
\end{cor}
\begin{proof}
Let $\{Z_{1},\dots,Z_{n}\}$ be a frame satisfying 1-2. Putting
$Z_{k}=X_{k}-iJX_{k}$ and $Z_{i}=X_{i}-iJX_{i}$, by hypothesis 1 it follows
$$
0=\n_{X_i-iJX_i}(X_k+iJX_k)(o)=\n_{X_{i}}X_k(o)+
\n_{JX_i}JX_k(o)+i(\n_{X_i}JX_k(o)-\n_{JX_{i}}X_{k}(o))\,.
$$
Then hypothesis 2 implies
\begin{equation}
0=\n_{JX_k}X_i(o)=\n_{X_k}JX_i(o)=\n_{JX_k}JX_i(o)\,.
\end{equation}
Finally we have
$$
\n_{Z_i}Z_k(o)=\n_{X_i}X_{k}(o)-\n_{JX_i}JX_k(o)
+i(\n_{X_{i}}JX_k(o)+\n_{JX_i}X_k(o))
$$
and so  $\n_{Z_i}Z_k(o)=0$ which implies $N_J=0$ .
\end{proof}
\section{Construction of generalized normal holomorphic frames in symplectic
  manifolds}
From \cite{J} we have the following
\begin{lemma}[Special frame on symplectic manifolds]
Let $(M,\kappa)$ be a 2$n$-manifold equipped with a
non-degenerate 2-form, $J\in\mathcal{C}_{\kappa}(M)$ and $g=g_{J}$.\\
The following facts are equivalent:
\begin{enumerate}
\item [1.]$\partial_{J}\ka=\overline{\partial}_{J}\ka=0$.
\vspace{0.1 cm}
\item [2.]For every $o \in M$ there exists a  local complex (1,0)-frame
      $\{Z_1,\dots, Z_n \}$ around $o$ such that:
      \begin{enumerate}
      \item [a)]$[Z_i ,Z_j ](o)=-(1/4)N_{J}[o](Z_i ,Z_j )$,
      $[\overline{Z}_{i},{Z}_{j}](o)=0$,\\ $1\leq r,s \leq n$;
      \vspace{0.1 cm}
      \item[b)] if
      $G_{ik}:=(g(Z_{i},\overline{Z}_{k}))$, then
      $G_{ik}(o)=\delta_{ik}$,  $d G_{ik}[o]=0$.
      \end{enumerate}
\vspace{0.1 cm}
\item[3.] For every $o\in M$, there exists a local complex (1,0)-coframe
      $\{\zeta_{1},\dots,\zeta_{n}\}$ around $o$ such that:
      \begin{enumerate}
      \item[a)] $\partial_{J}\zeta_{r}[o]=0$,
      $\overline{\partial}_{J}\zeta_{r}[o]=0$,
      $1\leq r \leq n$;
      \vspace{0.1 cm}
      \item[b)] if $H_{rs}:=(g(\zeta_{r},\overline{\zeta}_{s}))$,
      then $H_{rs}(o)=\delta_{rs}$, $d H_{rs}[o]=0$.
      \end{enumerate}
\end{enumerate}
\end{lemma}
We call a frame which satisfies $a)$ and $b)$ a \emph{``special frame''}.\\
Now we have
\begin{theorem}[Generalized normal holomorphic frames]
Let $(M,\ka)$ be a symplectic manifold, $J\in\mathcal{C}_{\kappa}(M)$,
$g=g_J$ and $o$ an arbitrary point of $M$.\\
Then there exists always a local complex (1,0)-frame,
$\{W_{1},\dots ,W_{n}\}$ around $o$, such that:
\begin{enumerate}
\item [1.]$\nabla_{W_{k}}\overline{W}_{i}(o)=0$, $1\leq k,i \leq n$.
\vspace{0.1 cm}
\item [2.]$\nabla_{W_{k}}W_{i}(o)$ $\in TM^{(0,1)}$, $1\leq k,i \leq n$.
\vspace{0.1 cm}
\item [3.]If $\;G_{rs}:=g(W_{r},\overline{W}_{s})$, then:
      $G_{rs}(o)=\delta_{rs}$, $d G_{rs}[o]=0$.
\vspace{0.1 cm}
\item [4.]$\nabla_{W_{r}}\nabla_{\overline{W}_{k}}W_{i}(o)=0$,
$1\leq r,k,i\leq n$.
\end{enumerate}
\end{theorem}
We call $\{W_{1},\dots ,W_{n}\}$ a \emph{generalized normal holomorphic frame}.
\begin{proof}
Step 1:\\
Let $\{Z_1,\dots,Z_n\}$ be a special frame around $o$.
\begin{enumerate}
\item[1.] From the relation
$$
0=[Z_{r},\ov{Z}_{k}](o)=
\n_{Z_{r}}\ov{Z}_{k}(o)-\n_{\ov{Z}_{k}}Z_{r}(o)
$$
it follows
$$
\n_{Z_{r}}\ov{Z}_{k}(o)=\n_{\ov{Z}_{k}}Z_{r}(o)
$$
and so by corollary 1
\begin{equation}
\n_{\ov{Z}_{k}}Z_{r}(o)=0,\qquad\n_{Z_{r}}\ov{Z}_{k}(o)=0
\quad 1\leq r,k\leq n\,.
\end{equation}
\item[2.] From (5), we have
\begin{equation}
\n_{Z_k}Z_i(o)\in T_oM^{(0,1)}\qquad 1\leq k,i\leq n\,.
\end{equation}
\item [3.]From (5) and corollary 1 we get that:
$$
g(\nabla_{Z_{r}}\nabla_{\ov{Z}_{k}}Z_i,Z_{s})(o)=
-g(\nabla_{\ov{Z}_{k}}Z_{i},\nabla_{Z_{r}}Z_{s})(o)+
Z_{r}g(\nabla_{\ov{Z}_{k}}Z{i},Z_{s})(o)=0\,;
$$
therefore
$\nabla_{Z_{r}}\nabla_{\ov{Z}_{k}}Z_{i}(o)$ are vector fields of type
(1,0) for any  $1\leq r,k,i\leq n$.
\item[4.] In a similar way:
$$
g(\nabla_{Z_{r}}\nabla_{\ov{Z}_{k}}Z_{i},\ov{Z}_{s})(o)=
Z_{r}g(\nabla_{\ov{Z}_{k}}Z_{i},\ov{Z}_{s})(o)\,.
$$
\end{enumerate}
So we can reduce to find a special frame $\{W_1,\dots,W_n\}$ around $o$  which
satisfies the relation
$$
\partial_{J}(g(\nabla_{\ov{W}_{k}}W_{i},\ov{W}_{s}))[o]=0\,.
$$
Step 2:\\
We can assume that a special frame satisfies at $o$
$$
Z_{i}=\de{}{z_{i}}\,,
$$
for some complex local coordinates $\{z_{1},\dots,z_{n}\}$ such that
$$
z_i(o)=0,\;\; 1\leq i \leq n\,.
$$
Let $\{W_1,,\dots,W_n\}$ be the complex frame
\begin{equation}
W_{i}=Z_{i}-\sum_{h=1}^{n}A_{hi}Z_{h}\,,
\end{equation}
where
$$
A_{si}=\sum_{a,b=1}^{n}Z_{a}\Gamma_{\ov{b}i}^{s}(o)z_{a}\ov{z}_{b}
$$
and
$$
\Gamma_{\ov{b}i}^{s}:=g(\nabla_{\ov{Z}_{b}}Z_{i},\ov{Z}_{s})\,.
$$
 It's easy
to check that this frame is a special one and then, by the first step of the
proof, it satisfies 1-3. So it's enough to show that $\{W_1,\dots,W_n\}$
satisfies 4.\\
Set $B_{si}=\delta_{si}-A_{si}$. Then we have
\begin{enumerate}
\item[1.] $B_{si}(o)=\delta_{si}$ ,
\item[2.] $\frac{\partial}{\partial z_a} B_{si}[o]=
       \frac{\partial}{\partial \ov{z}_b}B_{si}[o]=0$ ,
\end{enumerate}
then
$$
g(\nabla_{\ov{W}_{k}}W_{i},\ov{W}_{s})
=\sum_{a,b,c=1}^{n}\:\ov{B}_{ak}\ov{B}_{cs}
\ov{Z}_{a}(B_{bi})g(Z_{b},\ov{Z}_{c})+\ov{B}_{ak}\ov{B}_{cs}B_{bi}
g(\nabla_{\ov{Z}_{a}}Z_{b},\ov{Z}_{c})
$$
and
$$
\begin{aligned}
\partial_{J}(g(\nabla_{\ov{W}_{k}}W_{i},\ov{W}_{s}))[o]=&
\partial_{J}(\Gamma_{\ov{k}i}^{s})[o]+
\sum_{a=1}^{n}\ov{Z}_{a}(B_{si})\partial_{J}(\ov{B}_{ak})[o]+\\
&\sum_{b=1}^{n}\ov{Z}_{k}(B_{si}) \partial_{J}(\ov{B}_{bs})[o]+
\partial_{J}(\ov{Z}_{k}(B_{si}))[o]\,.
\end{aligned}
$$
Therefore
$$
\partial _{J}(g(\nabla_{\ov{W}_{k}}W_{i},\ov{W}_{s}))[o]=
\partial _{J}(\ov{Z}_{k}(B_{si}))[o]+\partial_{J}(\Gamma_{\ov{k}i}^{s})[o]\,.
$$
Hence we have
$$
\begin{aligned}
-\partial_{J}(\ov{Z}_{k}(B_{si}))[o]=&
\sum_{a,b=1}^{n}\partial_{J}(\ov{Z}_{k}(Z_{a}\Gamma_{\ov{b}i}^{s}(o)
z_{a}\ov{z}_{b}))[o]\\
&+\sum_{a,b=1}^{n}Z_{a}\Gamma_{\ov{b}i}^{s}(o)
\partial_{J}(\ov{Z}_{k}(z_{a}\ov{z}_b))[o]=\\
=&\sum_{a,b=1}^{n}Z_{a}\Gamma_{\ov{b}i}^{s}(o)
\partial_{J}(\ov{Z}_{k}(z_{a})\ov{z}_{b}+
\ov{Z}_{k}(\ov{z}_{b})z_{a})[o]\,.
\end{aligned}
$$
Finally we have
$$
Z_l(\ov{Z_{k}}(z_{a})\ov{z_{b}}+
\ov{Z_{k}}(\ov{z}_{b})z_{a})[o]=\delta_{kb}\delta_{la}\,.
$$
Therefore we obtain
$$
\partial_{J}(\ov{Z}_{k}(B_{si}))[o]=-\partial_{J}(\Gamma_{\ov{k}i}^{s})[o]
$$
which concludes the proof.
\end{proof}
\textbf{Remark:} In the previous lemma we can require
$$
\nabla_{W_k}\n_{W_j}\ov{W}_i(o)=0
$$
instead of
$$
\nabla_{W_{k}}\nabla_{\ov{W}_{j}}W_i(o)=0\,,
$$
but in general we can't require that the two conditions hold
simultaneously.\\

\textbf{Remark:} If $M$ is K\"ahler, then  the vector fields $\n_{W_i}W_j$
are globally of (1,0)-type and then they vanish at $o$.
Therefore in the K\"ahler case generalized normal holomorphic frames
are normal holomorphic frames.
\section{Integration of complex structures calibrated by symplectic forms}
In this section we apply the construction of the generalized normal
holomorphic frames (given in section 3) in order to find integrability
conditions for calibrated complex structures.\\
Let $(M,g,J)$ be a Hermitian manifold;
let us denote by $B$ the (2,1)-tensor defined by
$$
B(X,Y)=J(\nabla_{X}J)Y-(\nabla_{JX}J)Y\,.
$$
It is easy to check that $B$ satisfies
\begin{itemize}
\item[1.] $B(X,Y)-B(Y,X)=-N_J(X,Y)$;
\vspace{0.1 cm}
\item[2.] $B(Z_1,Z_2)=2iJ\nabla_{Z_1}Z_{2}+2\nabla_{Z_1}Z_2$;
\vspace{0.1 cm}
\item[3.] $B(Z_1,Z_2)\in TM^{(0,1)}$;
\vspace{0.1 cm}
\item[4.] $B(Z_1,\overline{Z}_2)=0$ ,
\end{itemize}
for every $ X,Y$ in $TM$ and $Z_1,Z_2$ in $TM^{(1,0)}$.\\
Therefore, we have:
\begin{lemma}
Let $(M,g,J)$ be a Hermitian manifold, then:
$$
N_J=0\iff B=0\,.
$$
\end{lemma}
In the symplectic case, it is natural to ask if the condition $B=0$  can be
weaken. An answer to this question is given by the following
\begin{theorem}
Let $(M,\kappa)$ be a symplectic manifold,
$J\in\mathcal{C}_{\kappa}(M)$ and $g=g_J$.
$$
\mbox {If }\; \nabla ''B=0 \mbox{ then M is a K\"ahler manifold},
$$
where $\n''$ is the (0,1)-part of the Levi-Civita connection
(i.e $\n''_Z:=\n_{Z^{(0,1)}}$).
\end{theorem}
\begin{proof}
Let $o$ be an arbitrary point of $M$ and let  $\{Z_1,\dots,Z_n\}$ be a
generalized  normal
holomorphic frame around $o$.\\
We have
\begin{equation*}
B(Z_i,Z_k)(o)=4\nabla_{Z_i}Z_k(o),\;\; 1\leq i,k\leq n\,.
\end{equation*}
By the properties of $B$ we obtain
\begin{equation*}
Z_l g(B(Z_i,Z_k),\ov{Z}_r)[o]=0,\;\; 1\leq l,i,k,r\leq n\,.
\end{equation*}
Then
\begin{equation*}
\begin{aligned}
0=&\ov{Z}_lg(B(Z_i,Z_k),\ov{Z}_r)[o]=\\
=&g(\nabla_{\ov{Z}_l}(B(Z_i,Z_k)),\ov{Z}_r)(o)+
g(B(Z_i,Z_k),\nabla_{\ov{Z}_l}\ov{Z}_r)(o)=\\
=&g((\nabla_{\ov{Z}_l}B)(Z_i,Z_k),\ov{Z}_r)(o)+
g(B(Z_i,Z_k),\nabla_{\ov{Z}_l}\ov{Z}_r)(o)=\\
=&g(B(Z_i,Z_k),\nabla_{\ov{Z}_l}\ov{Z}_r)(o)=\\
=&4g(\nabla_{Z_i}Z_k,\nabla_{\ov{Z}_l}\ov{Z}_r)(o)\,,
\end{aligned}
\end{equation*}
so in particular we obtain
\begin{equation*}
g(\nabla_{Z_i}Z_k,\nabla_{\ov{Z}_i}\ov{Z}_k)(o)=0\,.
\end{equation*}
Hence we get $\nabla_{Z_i}Z_k(o)=0$ for $1\leq i,k\leq n$ ,
which implies $N_J=0$.
\end{proof}
\textbf{Remark:} In  the previous theorem it's enough to
require
$$
g((\n_{\ov{Z}_1}B)(Z_1,Z_2),\ov{Z}_2)=0
$$
for every $Z_1,Z_2$ (1,0)-fields.\\

Now we give an integrability condition in terms of curvature.
\begin{lemma}
Let $(M,\kappa)$ be a  symplectic manifold,
 $J\in\mathcal{C}_{\kappa}(M)$,
A$g=g_J$. Let $R$ be  the curvature of $g$.
These facts are equivalent:
\begin{enumerate}
\item[1.] $\n''B=0\,;$
\vspace{0.1 cm}
\item[2.] $
R(\ov{Z},W)H=-i\n_{\ov{Z}}J\n_W H+
       iJ\n_W\n_{\ov{Z}}H+iJ\n_{\n_{\ov{Z}}W}H
       +\nabla_{\nabla_{W}\ov{Z}}H\,,
       \forall Z,W,H\in TM^{(1,0)}\,;
$
\vspace{0.1 cm}
\item[3.] $(\n_{\ov{Z}}J)(\n_{W}H)=0
    \quad\forall Z,W,H,S\in TM^{(1,0)}\,.$
\end{enumerate}
\end{lemma}
\begin{proof}
1$\iff$2:
Let $Z,W,H\in TM^{(1,0)}$, then
$$
\begin{aligned}
\frac{1}{2}(\nabla_{\ov{Z}}B)(W,H)=&\frac{1}{2}\{\nabla_{\ov{Z}}(B(W,H))
-B(\nabla_{\ov{Z}}W,H)-B(W,\nabla_{\ov{Z}}H)\}=\\
=&\n_{\ov{Z}}\n_{W}H+i\n_{\ov{Z}}J\n_{W}H-\n_{\n_{\ov{Z}}W}H-iJ
\n_{\n_{\ov{Z}}W}H\\
&-\n_{W}\n_{\ov{Z}}H-iJ\n_{W}\n_{\ov{Z}}H=\\
=&[\n_{\ov{Z}},\n_{W}]H+i\n_{\ov{Z}}J\n_{W}H-\n_{\n_{\ov{Z}}W}H-iJ
\n_{W}\n_{\ov{Z}}H\\
&-iJ\n_{\n_{\ov{Z}}W}H=\\
=&R(\ov{Z},W)H+i\n_{\ov{Z}}J\n_W H-
       iJ\n_W\n_{\ov{Z}}H\\
&-iJ\n_{\n_{\ov{Z}}W}H-\n_{\n_{W}\ov{Z}}H\,,
\end{aligned}
$$
and so
$$
\begin{aligned}
\n''B=0\iff  R(\ov{Z},W)H=&-i\n_{\ov{Z}}J\n_W H+iJ\n_W\n_{\ov{Z}}H\\
                          &+iJ\n_{\n_{\ov{Z}}W}H+\nabla_{\nabla_{W}\ov{Z}}H
\end{aligned}
$$
\hspace{3.0cm}$\forall Z,W,H\in TM^{(1,0)}$.\\
\newline
1$\iff$3: Let $o$ be a point in $M$ and let  $\{Z_1,\dots,Z_n\}$ be a
generalized
normal holomorphic frame around $o$. At the point $o$ we have
$$
\begin{aligned}
g((\n_{\ov{Z}_i}J)\n_{Z_j}Z_k,\ov{Z}_r)=&-g(J\n_{\ov{Z}_i}\n_{Z_j}Z_k,\ov{Z}_r)
+
g(\n_{\ov{Z}_i}J\n_{Z_j}Z_k,\ov{Z_r})=\\
&-ig(\n_{\ov{Z}_i}\n_{Z_j}Z_k,\ov{Z}_r)+i\ov{Z}_ig(\n_{Z_j}Z_k,\ov{Z}_r)\\
&+ig(\n_{Z_j}Z_k,\n_{\ov{Z}_i}\ov{Z}_r)=\\
=&2ig(\n_{Z_j}Z_k,\n_{\ov{Z}_i}\ov{Z}_r)
\end{aligned}
$$
for every  $1\leq i,j,k,r\leq n$. It is enough to prove 1$\iff$3.
\end{proof}
Recall that an Hermitian manifold $(M,g,J)$ is K\"ahler if and only
if $J$ is parallel, i.e.
$$
\n J=0\,.
$$
As application of Lemma 3 we have the following
\begin{theorem}
Let  $(M,\ka)$ be a symplectic manifold, $J\in\mathcal{C}_{\kappa}(M)$
and $g=g_J$. If
$$
(\n J)\n=0\,,
$$
then $(g,J,\ka)$ is a K\"ahler structure on $M$.
\end{theorem}
\vspace{1 cm}
\textbf{Acknowledgments:} The author is grateful to Paolo de
Bartolomeis who has proposed him the problems faced in this paper.\\
He is also grateful to Adriano Tomassini for his help in composing
this paper.

\end{document}